\documentclass[journal]{IEEEtran}    
\newcommand*\rot{\rotatebox{90}}
\def\smod{\color{black}}
\def\emod{\color{black}}
\usepackage[gen]{eurosym}
\usepackage{cite}
\usepackage{dsfont}

\usepackage{multicol}
\usepackage{balance}
\usepackage{caption}
\usepackage{mathtools}
\usepackage{mathrsfs}
\usepackage{subfigure}
\usepackage{graphicx,color}
\usepackage{latexsym}
\usepackage{amsmath} 
\usepackage{amssymb}  
\usepackage{amsmath,amssymb,amsfonts,theorem}
\usepackage{pgf}
\usepackage{tikz}
\usetikzlibrary{arrows,shadows,calc,decorations.shapes,decorations} 
\usepackage[underline=true,rounded corners=false]{pgf-umlsd}
\newtheorem{remark}{Remark}
\newtheorem{definition}{Definition}

\newtheorem{theorem}{Theorem}
\newtheorem{lemma}{Lemma}

\newtheorem{assumption}{Assumption}
\newtheorem{corollary}{Corollary}
\def\be{\begin{align}}
\def\ee{\end{align}}

\def\ba{\begin{split}}
\def\ea{\end{split}}

\title{
Optimal frequency regulation in nonlinear power networks including turbine-governor dynamics}
\author{Sebastian Trip$^{1}$  and Claudio De Persis$^{1}$
\thanks{$\star${The work of Sebastian Trip and Claudio De Persis is supported by the Danish Council for Strategic Research (contract no. 11-116843) within the `Programme Sustainable Energy and Environment', under the “EDGE” (Efficient Distribution of Green Energy) research project. The work of Claudio De Persis is also supported by the NWO (Netherlands Organisation for Scientific Research) programme Uncertainty Reduction in Smart Energy Systems (URSES) under the auspices of the project ENBARK. The results have appeared in \cite{trip_2017_tcns}. Preliminary results have appeared in \cite{trip_2016_acc} and \cite{trip_2016_ecc}. 
}
}
\thanks{$^{1}$Sebastian Trip and Claudio De Persis are with ENTEG, Faculty of Science and Engineering, University of Groningen, Nijenborgh 4, 9747 AG Groningen, the Netherlands. {\tt\small \{s.trip, c.de.persis\}@rug.nl}.}
}

\begin{document}
\maketitle
\thispagestyle{empty}
\pagestyle{empty}

\begin{abstract}
Motivated by an increase of renewable energy sources we propose a distributed \emph{optimal} Load Frequency Control scheme achieving frequency regulation and economic dispatch. Based on an energy function of the power network we derive an incremental passivity property for a well known nonlinear structure preserving network model, differentiating between generator and load buses. Exploiting this property we design \emph{distributed} controllers that adjust the power generation.  Notably, we explicitly include the turbine-governor dynamics where first-order and the widely used second-order dynamics are analyzed in a unifying way. Due to the non-passive nature of the second-order turbine-governor dynamics, incorporating them is challenging and we develop a suitable dissipation inequality for the interconnected generator and turbine-governor. This allows us to include the generator side more realistically in the stability analysis of optimal Load Frequency Control than was previously possible.
%
\end{abstract}

\section{Introduction}
Whenever there is an imbalance between generation and load, the frequency in the power network deviates from its nominal value. This makes frequency regulation, or `Load Frequency Control' (LFC), a critical task to maintain the stability of the network. Whereas primary droop control is utilized to act fast on smaller fluctuations to prevent destabilization, the frequency in the power network is conventionally regulated by `Automatic Generation Control' (AGC) that acts on the reference setting of the governors. To do so, each control area determines its `Area Control Error' (ACE) and changes the setpoints accordingly to compensate for local load changes and to maintain the scheduled tie-line power flows between different areas \cite{machowski_power_2008}, \cite{wood1996power}. However, due to an ever increasing penetration of renewable energy it is uncertain if the current AGC implementations are still adequate \cite{apostolopoulou_2016_tps}. The use of smart grids, computer-based control and communication networks offer on the other hand possibilities to improve the current practices \cite{pandey_2013_rser}, \cite{ibraheem_2005_tps}. Various solutions have been proposed to improve the performance of the AGC \cite{pan_1989_tps, liu_2012_pesgm, zhang_2013_tps, mi_2013_tps, yousef_2014_tps}. Specifically the effect of a large share of volatile renewable energy sources has been investigated \cite{variani_2013_tps}, \cite{xu_2016_tps}. Economic efficiency over slower timescales is  achieved by a tertiary optimization layer, commonly called the economic dispatch, that is outside of the conventional LFC loop.
\bigskip \\
Since the AGC was designed to be completely decentralized where each control area only reacts to its own ACE, there is loss of economic efficiency on the fast timescales of LFC. Instead of enforcing a predefined power flow over tie-lines, it is cost effective to coordinate the various regulation units within the whole system. This becomes especially relevant with a larger share of renewable energy sources where generation cannot be as accurately predicted as in the past. It is therefore desirable to further merge the secondary LFC and the tertiary optimization layer, which we call `\emph{optimal} Load Frequency Control' (OLFC). Although some centralized control schemes have been proposed \cite{DORFLER2017296,xi2017power}, the majority of current research focusses on distributed control architectures. The proposed distributed solutions to obtain OLFC can be roughly divided into two approaches. The first approach formulates the Lagrangian dual of the economic dispatch problem and solves the optimization problem based on a distributed primal-dual gradient algorithm that runs in parallel with the network dynamics \cite{zhang_2015_automatica, chen2015, stegink_2016_arxiv, you_2014_cdc, kasis_2016_arxiv, jokic_2009_epes, mudumbai_2012_ps, miao_2016_tps, cai_2015_cdc, apostolopoulou_2015_naps,Yi201545, cherukuri2016initialization, Yi2016259}. The advantage of this approach is that capacity constraints and convex cost functions can be straightforwardly incorporated. A drawback is however that generally information on the amount of uncontrollable generation and load needs to be available, which is generally unknown in LFC where only the frequency is used as a proxy for the imbalance. This issue is alleviated by the second approach, realizing that in the unconstrained case the marginal costs of the various generation units are identical at a cost effective coordination. In this approach optimality is achieved by employing a distributed consensus algorithm that converges to a state of identical marginal costs \cite{burger.et.al.mtns14b,trip_2016_automatica,schiffer_2016_ecc,zhao_2015_acc,monshizadeh_2015b_arxiv,andreasson_2013_ecc,kar_2012_pesgm,binetti_2014_tps,rahbari_2014_tsg,yang_2013_tps, yang2016, zhang_2012}.
 Although OLFC has been proposed as viable alternative to the conventional AGC, it poses the fundamental question if incorporating the economic dispatch into the LFC deteriorates the stability of the power network \cite{alvarado_2001_tps}.
\bigskip \\
\emph{Main Contributions.} This work continues and extends the study of the closed loop stability of OLFC and the power network. Specifically on the generation side there are still remaining challenges to include realistic models required in the study of frequency regulation. Recent advances in the analysis of OLFC in closed loop with the power network enable stability studies in presence of detailed generator models \cite{trip_2016_automatica}, \cite{stegink_2016B_arxiv} and improved network representations \cite{schiffer_2016_ecc}. However, including the important turbine-governor dynamics is less understood. We notice that indeed all of the referred studies on AGC \cite{pan_1989_tps, liu_2012_pesgm, zhang_2013_tps, mi_2013_tps, yousef_2014_tps, variani_2013_tps, xu_2016_tps} include a second-order model for the turbine-governor dynamics, whereas none of analytical studies on the stability of OLFC include such dynamics and are generally restricted to at most a first-order model. This paper makes the noteworthy extension towards closing this gap and incorporates the second-order turbine-governor dynamics in the stability analysis of the OLFC. We do this by establishing an incremental passivity property \cite{burger.et.al.mtns14b}, \cite{trip_2016_automatica} for a well studied structure-preserving network that represents various relevant power network configurations  \cite{bergen_1981_tpas}. This crucial passivity property of the power network is then exploited to incorporate first-order and second-order turbine-governor models in a unifying way. Including the second-order turbine-governor dynamics is especially challenging as they are non-passive and we cannot rely on the standard methodology for interconnecting passive systems. Instead, we develop a suitable dissipation inequality for the interconnected generator and turbine-governor. Due to the advantage of reduced generation and demand information requirements, we focus in this work on a distributed consensus based controller, where information on marginal costs is exchanged among neighbouring buses. Nevertheless we provide some guidelines how the higher order turbine-governor dynamics can be included in primal-dual based approaches as well. Along the stability analysis for the second-order turbine-governor model we establish a locally verifiable range of acceptable droop constants that allows us to infer frequency regulation. A case study confirms that a disregarding this range of droop constants in the controller design can lead to instability. We therefore argue that the design of an OLFC algorithm needs to carefully incorporate the effect of the turbine-governor dynamics. As a result of the distributed and modular design of the controllers, the proposed solution permits to straightforwardly include load control along the generation control and we provide a brief discussion on this topic.
\bigskip\\
The remainder of this paper is organized as follows.
In Section \ref{sec2}, we introduce the dynamic model of the power network, that we will study throughout this work.
In Section \ref{sec3}, we discuss the steady state of the power network and introduce an optimality criterium.
In Section \ref{sec4}, we prove an incremental cyclo-passivity property of the power network that is essential to the controller design.
In Section \ref{sec5}, we introduce the turbine-governor dynamics and propose distributed controllers that ensure frequency regulation and achieve economic dispatch.
In Section \ref{sec7}, we test our controllers in an academic case study using simulations.
In Section \ref{sec8}, conclusions and directions for future research are given.
\section{Power network model}\label{sec2}
We consider the nonlinear structure-preserving model of the power network proposed in \cite{bergen_1981_tpas} that we will extend in the later sections to include turbine-governor and load dynamics.
The network consists of $n_g$ generator buses and $n_l$ load buses. Each bus is assumed
 to be either a generator or a load bus,
such that the total number of buses in the network is $n_g + n_l = n$.
The network is represented by a connected and undirected graph
 $\mathscr{G} = (\mathcal{V}_g \cup \mathcal{V}_l, \mathcal{E})$, where $\mathcal{V}_g = \{1, \hdots, n_g\}$ is the set of generator buses, $\mathcal{V}_l = \{n_g +1, \hdots, n\}$ is the set of load buses and $\mathcal{E} = \{1,\hdots, m\}$ is
the set of transmission lines connecting the buses.
The network structure can be represented by its corresponding incidence matrix $\mathcal{B} \in \mathds{R}^{n \times m}$. The ends of transmission line $k$ are arbitrarily labeled with a `$+$' and a `$-$'. The incidence matrix is then given by
\[ \mathcal{B}_{ik} = \left\{
  \begin{array}{l l}
    +1 & \quad \text{if $i$ is the positive end of $k$}\\
    -1 & \quad \text{if $i$ is the negative end of $k$}\\
    0 & \quad \text{otherwise.}
  \end{array} \right.\]
  Following \cite{bergen_1981_tpas}, generator bus $i \in \mathcal{V}_g$ is modelled as
\begin{align}
\begin{split}
 \label{singlegenerator}
  \dot{\delta}_{i}=&~ \omega_{gi} \\
  M_{i}\dot{\omega}_{gi} =& -D_{gi}\omega_{gi} \\ &-  \sum_{j \in \mathcal{N}_{i}} V_{i}V_{j}B_{ij} \sin(\delta_{i} - \delta_{j}) + P_{mi},
\end{split}
\end{align}
  where $\mathcal{N}_i$ is the set of buses connected to bus $i$. \smod In high voltage tranmission networks considered here, the conductance is close to zero and therefore neglected, i.e. we assume the network to be lossless. \emod
 The uncontrollable loads\footnote{Controllable loads can be incorporated as well. The discussion on this topic is postponed to Remark \ref{remarkload} to facilitate a concise treatment.} are assumed \cite{bergen_1981_tpas}, \cite{kundur1994power} to consist of a constant and a frequency dependent component. We model a load bus for $i \in \mathcal{V}_{l}$ therefore as
\begin{align}
\begin{split}
  \dot{\delta}_{i} =&~\omega_{li}  \\
   0 =& -D_{li}\omega_{li} \\& -  \sum_{j \in \mathcal{N}_{i}} V_{i}V_{j}B_{ij} \sin(\delta_{i} - \delta_{j}) - P_{li}.
 \end{split}
\end{align}
  An overview of the used symbols is provided in Table 1.
 \begin{table}\label{tab1}
\setlength{\tabcolsep}{5pt}
\center
\medskip
\begin{tabular} {l  l}
&\bf State variables  \\
\hline\noalign{\smallskip}
  $\delta_{i}$ & Voltage angle \\
  $\omega_{gi}$ &  Frequency deviation at the generator bus \\
  $\omega_{li}$ & Frequency deviation at the load bus \\

  \noalign{\smallskip}
  &\bf Parameters \\
  \hline\noalign{\smallskip}
  $M_{i}$ & Moment of inertia \\
  $D_{gi}$ & Damping constant of the generator  \\
  $D_{li}$ & Damping constant of the load \\
  $B_{ij}$ & Susceptance of the transmission line\\
  $V_i$ & Voltage \\
\noalign{\smallskip}
&\bf Controllable input \\
\hline\noalign{\smallskip}
  $P_{mi}$   & Mechanical power \\
   &\bf Uncontrollable input \\
\hline\noalign{\smallskip}
      $P_{li}$ & Unknown constant power demand \\
\end{tabular}
\caption*{Table 1: Description the variables and parameters appearing in the power network model.}
\end{table}
Since the power flows are determined by the differences in voltage angles, it is convenient to introduce $\eta_k = \delta_i - \delta_j$, where $\eta_k$ is the difference of voltage angles across line $k$ joining buses $i$ and $j$.
 For all buses the dynamics of the power network are written as
  \begin{align}
\begin{split}\label{syscompact}
   \dot{\eta} =&~\mathcal{B}^T \omega \\
  M\dot{\omega}_{g}  =& -D_{g}\omega_{g} -  \mathcal{B}_{g} \Gamma \sin(\eta) + P_m\\
   \boldsymbol{0} =&  -D_{l}\omega_{l} -  \mathcal{B}_{l} \Gamma \sin(\eta) - P_l,\\
  \end{split}
\end{align}
  where   $\omega =  (\omega_{g}^T, \omega_{l}^T)^T$, $\eta= \mathcal{B}^T \delta$ and $\Gamma={\rm diag}\{\gamma_1,\ldots, \gamma_m\}$, with $\gamma_k=V_i V_j B_{ij} = V_j V_i B_{ji}$ and the index $k$ denoting the line $\{i,j\}$. The matrices $\mathcal{B}_{g} \in \mathds{R}^{n_{g} \times m}$ and $\mathcal{B}_{l} \in \mathds{R}^{n_{l} \times m}$ are obtained by collecting from $\mathcal{B}$ the rows indexed by $\mathcal{V}_{g}$ and $\mathcal{V}_{l}$ respectively. The remaining symbols follow straightforwardly from the node dynamics and are diagonal matrices or vectors of suitable dimensions. It is possible to eliminate $\omega_{l}$ in (\ref{syscompact}) by exploiting the identity $\omega_{l} = D_{l}^{-1}(-  \mathcal{B}_{l} \Gamma \sin( \eta) - P_l)$ and realizing that $\mathcal{B}^T \omega = \mathcal{B}^T_{g} \omega_{g} + \mathcal{B}^T_{l} \omega_{l}$ \cite{monshizadeh_2015_arxiv}. As a result we can write (\ref{syscompact}) equivalently as
\begin{align}
\begin{split} \label{syscompactred2}
   \dot{\eta} =&~\mathcal{B}_{g}^T \omega_{g} +\mathcal{B}_{l}^TD_{l}^{-1}(-  \mathcal{B}_{l} \Gamma \sin(\eta) - P_l ) \\
  M\dot{\omega}_{g}  =& -D_{g}\omega_{g} -  \mathcal{B}_{g} \Gamma \sin(\eta) + P_m.
\end{split}
\end{align}
  We will however keep $\omega_l$ when it enhances the readability of this paper.
  \begin{remark}[Control areas]
   In the absence of load buses, the considered model appears in the study of automatic generation control of control areas, where a control area is described by an equivalent generator. A control area is then typically modelled as
    \begin{align}
\begin{split}
  \dot{\delta}_{i} =&~\omega_{gi} \\
  M_{i}\dot{\omega}_{gi} =& -D_{gi}\omega_{gi} \\ &-  \sum_{j \in \mathcal{N}_{i}} V_{i}V_{j}B_{ij} \sin(\delta_{i} - \delta_{j}) + P_{mi} -P_{li},
 \end{split}
\end{align}
  where the loads are collocated at the equivalent generator.
  All results in this paper also hold for this particular case.
  \end{remark}
\begin{remark}[Microgrids]
 Besides modelling high voltage power networks, system (\ref{syscompact}) has also been used to model (Kron reduced)
microgrids \cite{schiffer_2015_phd, trip_2014_cdc, dorfler_2016_cns, shafiee_2015_pe, persis_2016_arxiv_microgrid, schiffer13_aut}. Smaller synchronous machines
and inverters are then represented by (1) and (2) respectively.
\end{remark}\smod

\begin{remark}[Detailed network models]
To stress the contribution of this work we focus on a basic structure preserving model of the power network. The voltages in this paper are considered constant, which is a common assumption in models tailored to study frequency regulation, since the voltage dynamics are (generally) fast compared to the frequency dynamics \cite{garcia_2012_book}, \cite{bevrani2009robust}. As becomes clear in the subsequent sections, our analysis depends mostly on the existence of an energy function for the considered model.
These energy functions have been developed for more realistic network models that e.g. include voltage dynamics, exciter dynamics and that distinguish between internal and terminal generator buses \cite{trip_2016_automatica}, \cite{stegink_2016B_arxiv}, \cite{chu_2005_cssp}, \cite{padiyar_2013_book}. Commonly these energy functions include a kinetic term $\frac{1}{2} \omega_g^T M \omega_g$, which in our work is essential to derive the passivity property that we exploit in the controller design. It is therefore expected that the proposed design can be extended to more complex network dynamics as well. Specifically, the passivity property derived for the model at hand (see Lemma 3) has explicitly been established for generators including voltage dynamics in \cite{trip_2016_automatica} and \cite{stegink_2016B_arxiv}.
\end{remark}
\emod
\section{Steady state and optimality} \label{sec3}
Before addressing the turbine-governor dynamics that adjust $P_m$, we discuss the steady state frequency deviation under constant generation $\overline P_m$. In particular we study the optimal value of $\overline P_m$ that allows for a zero frequency deviation at steady state, i.e. $\overline \omega  = \boldsymbol{0}$.
The steady state $(\overline \eta, \overline \omega, \overline P_m)$ of (\ref{syscompact}) necessarily satisfies
\begin{align}
\begin{split} \label{syscompactsteady}
   \boldsymbol{0} =&~\mathcal{B}^T \overline \omega \\
   \boldsymbol{0}  =& -D_{g}\overline \omega_{g} -  \mathcal{B}_{g} \Gamma \sin(\overline \eta) + \overline P_m \\
   \boldsymbol{0} =&  -D_{l}\overline \omega_{l} -  \mathcal{B}_{l} \Gamma \sin(\overline \eta) - P_l. \\
  \end{split}
\end{align}
  We make the natural assumption that a, \smod possibly non-unique, \emod solution to (\ref{syscompactsteady}) exists, which corresponds to the ability of the network to transfer the required power at steady state.
   \begin{assumption}[Solvability]\label{assum1}
    For a given $P_l \in \mathds{R}^{n_l}$ and $\overline P_m \in \mathds{R}^{n_{g}}$, there exist $\overline \eta \in {\rm Im}(\mathcal{B}^T)$, $\overline \omega \in {\rm Ker}(\mathcal{B}^T)$ such that (\ref{syscompactsteady}) is satisfied.
  \end{assumption}
  From algebraic manipulations of (\ref{syscompactsteady}) we can derive the following lemma that makes the frequency deviation at steady state $\overline \omega$ explicit.
 \begin{lemma}[Steady state frequency]\label{lemma1}
    Let Assumption \ref{assum1} hold, then necessarily $\overline \omega = \mathds{1}_n \omega_*$, with
    \begin{align}\label{omega.star}
 \omega_\ast =   \frac{ \mathds{1}_{n_g}^T \overline P_m- \mathds{1}_{n_l}^T P^{l} }{\mathds{1}_{n_{g}}^T D_{g} \mathds{1}_{n_{g}} + \mathds{1}_{n_{l}}^T D_{l} \mathds{1}_{n_{l}} }
,
\end{align} where $\mathds{1}_n \in \mathds{R}^{n}$ is the vector consisting of all ones.
  \end{lemma}
 We recover therefore the well known fact that the total generation needs to be equal to the total load in order to have a zero frequency deviation in a lossless network. As we only require the total generation to be equal to the total load, it is natural to wonder if we can distribute the generation in an optimal manner. To this end, we assign to every generator a strictly convex linear-quadratic cost function that relates the generated power $P_{mi}$ to the generation costs $C_i(P_{mi})$, typically expressed in \$/MWh, i.e.
\begin{align}
C_i(P_{mi}) = \frac{1}{2}q_{i}P_{mi}^2 + r_{i}P_{mi} + s_{i}.
\end{align}

To formalize the notion of optimality in this work, we pose the following optimization problem:
\begin{align}
\begin{split}\label{optimal}
& \min_{P_m} C(P_m)\\
{\rm s.t.} \quad & 0~= ~\mathds{1}_{n_g}^T \overline P_{m}- \mathds{1}_{n_l}^T P^{l},
\end{split}
\end{align}
where $C(P_m)=\sum_{i \in \mathcal{V}_{g}}C_i(P_{mi})$. Defining furthermore $Q = \text{diag}(q_1, \hdots, q_{n_g})$, $R = (r_1, \hdots, r_{n_g})^T$ and $S = (s_1, \hdots, s_{n_g})^T$ we can compactly write
\begin{align}C(P_m) =  \frac{1}{2}P_m^TQP_m + R^T P_m + \mathds{1}^T_{n_g}S.
\end{align}
From the discussion of Lemma 1, we note that satisfying the equality constraint in (\ref{optimal}) implies $\overline \omega = \boldsymbol{0}$. The solution to (\ref{optimal}), indicated by the superscript $opt$, therefore satisfies \cite[Lemma 4]{trip_2016_automatica}
\begin{align}
\begin{split}\label{optss}
   \boldsymbol{0} =&~ \mathcal{B}^T \boldsymbol{0}\\
  \boldsymbol{0}   =& - D_{g} \boldsymbol{0}-  \mathcal{B}_{g} \Gamma \sin(\overline \eta) + \overline P_m^{opt} \\
 \boldsymbol{0} =&  -D_{l}\boldsymbol{0} -  \mathcal{B}_{l} \Gamma \sin(\overline \eta) - P_l.\\
\end{split}
\end{align}
  It is possible to explicitly characterize the solution to (\ref{optimal}).
\begin{lemma}[Optimal generation]\label{lemma4}
The solution $\overline P^{opt}_m$ to (\ref{optimal}) satisfies
\begin{align}\label{optimal.u}
\overline P^{opt}_m =~ Q^{-1}(\overline \lambda^{opt} - R),
\end{align}
where
\begin{align} \label{optimal.u2}
\overline \lambda^{opt} =~ \frac{\mathds{1}_{n_g} (\mathds{1}_{n_l}^T P_l + \mathds{1}_{n_g}^TQ^{-1}R)}{\mathds{1}_{n_g}^T Q^{-1}\mathds{1}_{n_g}}. 
\end{align}
\end{lemma}

The first derivative of the cost function is commonly called the `marginal cost function'.
 From (\ref{optimal.u}) and (\ref{optimal.u2}) it is then immediate to see that
 \begin{align}\label{marginal}
 Q \overline P^{opt}_m+ R = \overline \lambda^{opt} \in {\rm Im}(\mathds{1}_{n_g}),
  \end{align}which implies that at the solution to (\ref{optimal}) all marginal costs are identical.
\begin{remark}[Information requirements]
 Solving (\ref{optimal}) explicitly requires the knowledge of the total load $\mathds{1}^T_{n_l}P_l$. A popular approach to solve (\ref{optimal}) in a distributed fashion is based on primal-dual gradient dynamics \cite{zhang_2015_automatica, chen2015, stegink_2016_arxiv, you_2014_cdc, kasis_2016_arxiv, jokic_2009_epes, mudumbai_2012_ps, miao_2016_tps, cai_2015_cdc, apostolopoulou_2015_naps}.  Commonly, these approaches do however require knowledge of the loads or power flows. A remarkable feature of our work is that the proposed distributed controllers, that will be discussed in the remaining of this paper, solve (\ref{optimal}) without such measurements at the cost of the restriction to linear-quadratic cost functions and the absence of generation and power flow constraints.
\end{remark}
The focus of this section was the characterization of the (optimal) steady state of the power network under constant power generation. In the next section we establish a passivity property of the power network that will be useful to design controllers that dynamically adjust $P_m$, ensuring that $P_m$ converges to the optimal steady state $\overline P_m^{opt}$.
\section{An incremental passivity property of the power network}\label{sec4}
We now establish a passivity property for the considered power network model, that is essential to the stability analysis in the following section. Being more specific, we show that (\ref{syscompact}) is output strictly incrementally cyclo-passive \cite{hill_1980_franklin}, \cite{willem_2007_ejc}, \cite{Pavlov2008}  with respect to its steady state solution, when we consider $P_m$ as the input and $\omega_g$ as the output. We first recall the following definition (with some abuse of terminology\footnote{\smod Incremental passivity as defined in e.g. \cite{Pavlov2008} holds for any two solutions to the system. In the definition here, incremental passivity is required to hold with respect to a steady state solution. \emod  }):
\begin{definition}[Incremental cyclo-passivity]
  System \begin{align} \begin{split}
  \dot x =&~ f(x,u)\\
   y=&~h(x),
  \end{split}
  \end{align}
  $x \in \mathcal{X}$, $\mathcal{X}$ the state space, $u,y \in \mathds{R}^n$, is incrementally cyclo-passive with respect to a constant triplet $(\overline x, \overline u, \overline y)$ satisfying
  \begin{align}
  \begin{split}
    \boldsymbol{0} =& ~f(\overline x, \overline u) \\
    \overline y =& ~h(\overline x),
    \end{split}
  \end{align}
  if there exists a continuously differentiable function $\mathcal{S} : \mathcal{X} \rightarrow \mathds{R}$, such that for all $x \in \mathcal{X}$, $u \in \mathds{R}^m$ and $y=h(x)$, $\overline y = h(\overline x)$
  \begin{align}
  \begin{split}
  \dot{\mathcal{S}} = \frac{\partial \mathcal{S}}{\partial x} f(x,u) + \frac{\partial \mathcal{S}}{\partial \overline x} f(\overline x,\overline u)  \leq& -\|y-\overline y\|^2_W \\&+(y-\overline y)^T (u-\overline u), \nonumber
  \end{split}
  \end{align}
  where $\|y-\overline y\|^2_W = (y - \overline y)^TW(y- \overline y)$.
  If $W > \boldsymbol{0}$, the system is output strictly incrementally cyclo-passive.
\end{definition}

We remark that the definition above differs from the ordinary definition of incremental passivity in that it includes the prefix `\emph{cyclo-}' indicating that $\mathcal{S}$ is not required to be positive definite nor to be bounded from below. If $\mathcal{S}$ is positive definite, we call the system incrementally passive. We now show that the power network satisfies Definition 1 above.
  \begin{lemma}[Incremental cyclo-passivity of (\ref{syscompactred2})]{\label{u1}}
Let Assumption 1 hold. System
(\ref{syscompactred2}) with input $P_m$ and output $\omega_g$ is an output strictly incrementally cyclo-passive system, with respect to $(\overline \eta, \overline \omega_g)$ satisfying
 \begin{align}
\begin{split}\label{syscompactred2ss}
   \boldsymbol{0} =&~\mathcal{B}_{g}^T \overline \omega_{g} +\mathcal{B}_{l}^TD_{l}^{-1}(-  \mathcal{B}_{l} \Gamma \sin(\overline \eta) - P_l ) \\
  \boldsymbol{0}  =& -D_{g}\overline \omega_{g} -  \mathcal{B}_{g} \Gamma \sin(\overline \eta) + \overline P_m. \\
 \end{split}
\end{align}
 Namely, there exists a storage function
$U(\eta, \overline \eta, \omega_{g}, \overline\omega_{g})$ which satisfies the following incremental dissipation inequality
\begin{align}
\begin{split}
\dot U  =& -\|\omega_{g} - \overline \omega_{g}\|_{ D_{g}}^2  -\|\omega_{l} - \overline \omega_{l}\|_{ D_{l}}^2 \\& + (\omega_{g} - \overline \omega_{g})^T(P_m - \overline P_m),
\end{split}
\end{align}

where $\dot U$ represents the derivative of $U(\eta, \overline \eta, \omega_{g}, \overline\omega_{g})$ along the solutions to  (\ref{syscompactred2}).
  \end{lemma}
  \begin{IEEEproof}
  Consider the incremental storage function
  \begin{align}\label{storagefunction}
\begin{split}
   U(\eta, \overline\eta, \omega_{g}, \overline \omega_{g}) =& ~\frac{1}{2}(\omega_{g} - \overline \omega_{g})^TM(\omega_{g} - \overline \omega_{g}) \\ & -\mathds{1}^T\Gamma\cos(\eta) + \mathds{1}^T \Gamma\cos(\overline \eta) \\&-  \left(\Gamma\sin(\overline \eta)\right)^T(\eta - \overline \eta).
  \end{split}
\end{align}

  We have that $U(\eta, \overline\eta,\omega_{g}, \overline \omega_{g})$ satisfies along the solutions to (\ref{syscompactred2})
    \begin{align}\label{derivativeU}
\begin{split}
  \dot U =&~ (\omega_{g} - \overline \omega_{g})^T (-D_{g}\omega_{g} -  \mathcal{B}_{g} \Gamma \sin(\eta) + P_m) \\
  &+(\Gamma\sin(\eta) - \Gamma\sin(\overline \eta))^T \\ &\hspace{1em}\cdot (\mathcal{B}_{g}^T  \omega_{g} +
  \mathcal{B}_{l}^TD_{l}^{-1}(-  \mathcal{B}_{l} \Gamma \sin(\eta) - P_l))
  \\
   =& -\|\omega_{g} - \overline \omega_{g}\|_{ D_{g}}^2
   + (\omega_{g} - \overline \omega_{g})^T(P_m - \overline P_m)
   \\ & + (\Gamma\sin(\eta) - \Gamma\sin(\overline \eta))^T\mathcal{B}_{l}^T(\omega_l - \overline \omega_l)
   \\ =& -\|\omega_{g} - \overline \omega_{g}\|_{ D_{g}}^2
   + (\omega_{g} - \overline \omega_{g})^T(P_m - \overline P_m)
   \\ & + (\Gamma\sin(\eta)\mathcal{B}_{l} + P_l )^T D_l^{-1} D_l(\omega_l - \overline \omega_l)
   \\ & - (\Gamma\sin(\overline \eta)\mathcal{B}_{l} +P_l)^T D_l^{-1} D_l(\omega_l - \overline \omega_l)
    \\ =&-\|\omega_{g} - \overline \omega_{g}\|_{ D_{g}}^2  -\|\omega_{l} - \overline \omega_{l}\|_{ D_{l}}^2 \\& + (\omega_{g} - \overline \omega_{g})^T(P_m - \overline P_m),
  \end{split}
\end{align}

  where we exploit identity (\ref{syscompactred2ss}) in the second equation.
  \end{IEEEproof}
  \bigskip
  Note that the result of Lemma \ref{u1} holds in particular if we take $\overline \omega = \boldsymbol{0}$ and $\overline P_m = \overline P_m^{opt}$.
  We now consider what conditions ensure that storage function (\ref{storagefunction}) has a local minimum at a steady state satisfying (\ref{syscompactred2ss}).
  \begin{assumption}[Steady state angle differences]\label{assum2}
The differences in voltage angles $\overline \eta$ in (\ref{syscompactsteady}) satisfy $\overline \eta_k \in (\frac{-\pi}{2},\frac{\pi}{2}) \quad \forall k \in \mathcal{E}$.
\end{assumption}
\smod Note that Assumption \ref{assum2} is  generally satisfied under normal operating conditions of the power network, where a small difference in voltage angle is also referred to as phase-cohesiveness \cite{dorfler2013} and is preferred to avoid instability after perturbations \cite{nerc2016}. \emod
 \begin{lemma}[Local minimum of (\ref{storagefunction})]\label{lemmabregman}
   Let Assumption \ref{assum2} hold. Then the storage function (\ref{storagefunction}) has a local minimum at $(\overline \eta, \overline \omega_g)$.
 \end{lemma}
 \begin{IEEEproof}
\smod
We first recall the definition of a Bregman distance \cite{bregman_1967_CMMP}. Let $F: \mathcal{X} \rightarrow \mathds{R}$ be a continuously differentiable and strictly convex function defined on a closed convex set $\mathcal{X}$. The Bregman distance associated with $F$ for the points $x, \overline x$ is defined as
\begin{align}
D_F(x,\overline x) = F(x) - F(\overline x) - \nabla F(\overline x)^T(x - \overline x).
\end{align}
A useful property of $D_F$ is that it is positive definite in its first argument, due to the strict convexity of $F$. Lemma \ref{lemmabregman} then follows from (\ref{storagefunction}) being the Bregman distance associated with the function $F(\eta, \omega_g) = \frac{1}{2}\omega_g^T M \omega_g - \mathds{1}_m^T \Gamma \boldsymbol{\cos}(\eta)$, which is strictly convex at the point $(\overline \eta, \overline \omega_g)$ under Assumption \ref{assum2}. \emod
 \end{IEEEproof}
\begin{remark}[Boundedness of solutions]
  In the proof of Theorem 1 we require Assumption \ref{assum2} and subsequently Lemma \ref{lemmabregman} to ensure that there exists a compact forward invariant set around an equilibrium of (\ref{syscompact}). This allows us to apply LaSalle's invariance principle in the stability analysis.
\end{remark}
In this section we have established that the power network model (\ref{syscompact}) is an output strictly incrementally cyclo-passive system. Furthermore we have shown that under Assumption 2, the incremental storage function $U$ has a local minimum at its steady state. These results turn out the be essential to the design of the distributed controllers in the next section and to prove asymptotic stability of the obtained closed-loop system.

\section{Optimal turbine-governor control}\label{sec5}
The generated power $P_{mi}$ at generator $i$ is the output of the turbine-governor system. Various turbine-goveror models appear in the literature. We consider two of the most widely used models that have fundamentally different properties. We therefore partition the set of generators $\mathcal{V}_{g} = \mathcal{V}_{g1} \cup \mathcal{V}_{g2}$ into the sets $\mathcal{V}_{g1}$ and $\mathcal{V}_{g2}$, where the turbine-governor dynamics are described by first-order and second-order dynamics respectively. Being able to incorporate both types in a single framework, unifies the various modelling assumptions appearing in conventional AGC and OLPC studies, and increases the modelling flexibility.
\medskip \\
The first-order and second-order turbine-governor dynamics will be discussed separately and controllers are proposed that achieve frequency regulation.
 To facilitate the controller design using only local information we write (\ref{derivativeU}), taking therein and in the remainder of this work $\overline \omega = \boldsymbol{0}$ and $\overline P_m = \overline P_m^{opt}$, as
\begin{align}
\begin{split}
  \dot U =&-\|\omega_{g} - \boldsymbol{0}\|_{ D_{g}}^2  -\|\omega_{l} -\boldsymbol{0}\|_{ D_{l}}^2 \\& + (\omega_{g} - \overline \omega_{g})^T(P_m - \overline P_m^{opt})\\
  =& \sum_{i \in \mathcal{V}_g} \dot U_{gi}(\omega_{gi},P_{mi},\overline P_{mi}^{opt})   + \sum_{i \in \mathcal{V}_l} \dot U_{li}(\omega_{li}),
  \end{split}
\end{align}
where we define with a slight abuse of notation
\begin{align}\label{ugi}
\begin{split}
\dot U_{gi}(\omega_{gi},P_{mi},\overline P_{mi}^{opt})  =& -D_{gi}\omega_{gi}^2 + \omega_{gi}(P_{mi} - \overline P_{mi}^{opt}) \\
 \dot U_{li}(\omega_{li}) =& -D_{li}\omega_{li}^2.
 \end{split}
 \end{align}
For the sake of exposition we only consider decentralized controllers in subsections \ref{subA} and \ref{subB} that guarantee frequency regulation
without achieving optimality.
These results are then instrumental to Subsection \ref{subC} where a distributed control architecture is proposed with  controllers that exchange information on their marginal costs with their neighbours over a communication network to achieve optimality.
\begin{table}\label{tab1}
\setlength{\tabcolsep}{5pt}
\center
\medskip
\begin{tabular} {l  l}
&\bf State variables  \\
\hline\noalign{\smallskip}
  $P_{si}$  & Steam power \\
  $P_{mi}$  & Mechanical power \\
  \noalign{\smallskip}
  &\bf Parameters \\
  \hline\noalign{\smallskip}
   $T_{si}$ & Governor time constant \\
   $T_{mi}$ & Turbine time constant \\
   $K_{i}$ & Droop constant \\
\noalign{\smallskip}
&\bf Controllable input \\
\hline\noalign{\smallskip}
  $\theta_{i}$ & Power generation control \\
\end{tabular}
\caption*{Table 2: Description of the variables and parameters appearing in the turbine-governor dynamics.}
\end{table}

\subsection{First-order turbine-governor dynamics}\label{subA}
We start with the first-order turbine-governor dynamics of a single generator $i \in \mathcal{V}_{g1}$. The
dynamics are given by
\begin{align}\label{1tg}
   T_{mi}\dot{P}_{mi}  =& -P_{mi} - K^{-1}_i\omega_{gi} + \theta_{i},
\end{align}
where $\theta_i$ is an additional control input to be designed. An overview of the used symbols is provided in Table 2.
Consider the following controller at bus $i$: \begin{align}\label{1tgc}
  \begin{split}
    T_{\theta_{i}} \dot \theta_{i} =& -\theta_{i} + P_{mi},
  \end{split}
  \end{align}
  where the controller time constant $T_{\theta_i}$ can be chosen to obtain a desirable rate of change of the control input $\theta_i$.
  As explained before, an additional communication term will be added to controller (\ref{1tgc}) in subsection \ref{subC}  to enforce optimality at steady state. The following lemma provides an intermediate result that is useful later on.
\begin{lemma}[Incremental passivity of (\ref{1tg}), (\ref{1tgc})]
System
(\ref{1tg}), (\ref{1tgc}) with input $-\omega_{gi}$ and output $P_{mi}$ is an incrementally passive system, with respect to $(\overline P_{mi}^{opt}, \overline \theta_i)$ satisfying
\begin{align}\ba\label{tb1ss}
   0  =& -\overline P^{opt}_{mi} - K^{-1}_i 0 + \overline \theta^{opt}_{i}\\
   0 =& -\overline \theta_{i}^{opt} + \overline P_{mi}^{opt},
\ea\end{align}
 Namely, there exists a positive definite storage function
$Z_{1i}(P_{mi}, \overline P_{mi}^{opt}, \theta_i, \overline \theta_i^{opt})$ which satisfies the following incremental dissipation inequality
\be
\ba \label{zi}
\dot Z_{1i} =& -K_i(\theta_i - P_{mi})^2 -\omega_{gi}(P_{mi}-\overline P_{mi}^{opt}),
\ea
\end{align}
where $\dot Z_{1i}$ represents the derivative of $Z_{1i}(\theta_i, \overline \theta_i^{opt}, P_{mi}, \overline P_{mi}^{opt})$ along the solutions to (\ref{1tg}), (\ref{1tgc}).
  \end{lemma}
  \begin{IEEEproof}
  Consider the incremental storage function
  \begin{align}
    Z_{1i}= \frac{T_{\theta_i}K_i}{2}(\theta_i - \overline \theta_i^{opt})^2 + \frac{T_{P_{mi}}K_i}{2}(P_{mi} - \overline P_{mi}^{opt})^2.
  \end{align} We note that $Z_{1i}$ satisfies along the solutions to (\ref{1tg}), (\ref{1tgc}),
  \begin{align}
  \begin{split}
  \dot Z_{1i} =& -(\theta_i - \overline \theta_i^{opt})K_i\theta_i + (\theta_i - \overline \theta_i^{opt})K_iP_{mi} \\
  & -(P_{mi} - \overline P_{mi}^{opt})K_iP_{mi} + (P_{mi} - \overline P_{mi}^{opt})K_i\theta_i\\& - (P_{mi} - \overline P_{mi}^{opt})\omega_{gi}\\
  =& -K_i(\theta_i - P_{mi})^2 - (P_{mi} - \overline P_{mi}^{opt})\omega_{gi},
  \end{split}
  \end{align}
  where we exploit identity (\ref{tb1ss}) in the second equation.
  \end{IEEEproof}
  \medskip
  The interconnection of generator dynamics (\ref{syscompact}) and turbine-governor dynamics (\ref{1tg}) including controller (\ref{1tgc}) can be understood as a feedback interconnection of two incrementally passive systems. The following corollary is then an immediate result from this observation.
  \begin{corollary}[Passive interconnection]
  Along the solutions to (\ref{syscompact}), (\ref{1tg}) and (\ref{1tgc}),  $Z_{1i}(\theta_i, \overline \theta_i^{opt}, P_{mi}, \overline P_{mi}^{opt})$ satisfies
  \begin{align}
  \begin{split}
    \dot U_{gi} + \dot Z_{1i} =& -D_{gi}\omega_{gi}^2  -K_i(\theta_i - P_{mi})^2 \leq 0,
    \end{split}
    \end{align}
    where $\dot U_{gi}$ and $\dot Z_{1i}$ are given in (\ref{ugi}) and (\ref{zi}) respectively.
  \end{corollary}
  We now perform a similar analysis for the second-order turbine-governor dynamics.
\subsection{Second-order turbine-governor dynamics}\label{subB}
Consider the second-order turbine-governor dynamics of a single generator $i \in \mathcal{V}_{g2}$. The
dynamics are given by
 \be\ba\label{2tg}
   T_{si}\dot{P}_{si}  =& -P_{si} - K^{-1}_i\omega_{gi}  + \theta_{i} \\
   T_{mi}\dot{P}_{mi}  =& -P_{mi} + P_{si}, \\
  \ea
  \end{align}
  where $\theta_i$ is again an additional control input to be designed. In contrast to the first-order dynamics, the second-order dynamics do not possess a useful passivity property. This can be readily concluded from the observation that system (\ref{2tg}) with input $\omega_{gi}$ and output $P_{mi}$ has relative degree 2. We now propose a different controller than (\ref{1tgc}) to accommodate the higher order turbine-governor model, namely
  \begin{align}\label{2tgc}
  \begin{split}
    T_{\theta_{i}} \dot \theta_{i} =& -\theta_{i} + P_{si} - (1 - K^{-1}_i)\omega_{gi},
  \end{split}
  \end{align}
  where $K^{-1}$ is the droop constant appearing in (\ref{2tg}).
  Similar to (\ref{1tgc}) we postpone adding an additional communication term until the next subsection.

  \begin{lemma}[Storage function for second-order dynamics]\label{lemma2ndorder}
   There exists a positive definite storage function
$Z_{2i}(P_{si}, \overline P_{si}, P_{mi}, \overline P_{mi}^{opt}, \theta_i, \overline \theta_i^{opt})$ which satisfies along the solutions to (\ref{syscompact}), (\ref{2tg}) and (\ref{2tgc})
  \begin{align}
  \begin{split}
    & \dot U_{gi} + \dot Z_{2i} =  \begin{bmatrix}
      \omega_{gi} \\ P_{si} - P_{mi} \\ P_{si} - \theta_i
    \end{bmatrix}^T  W_i \begin{bmatrix}
      \omega_{gi} \\ P_{si} - P_{mi} \\ P_{si} - \theta_i
    \end{bmatrix},
    \end{split}
    \end{align}
    with
    \be
    W_i = \begin{bmatrix}
      -D_{gi} & -\frac{1}{2}K_i^{-1} -\frac{1}{2} & -\frac{1}{2}K_i^{-1} + \frac{1}{2} \\ -\frac{1}{2}K_i^{-1} -\frac{1}{2}  &-T_{si}T_{mi}^{-1} &- \frac{1}{2} \\ -\frac{1}{2}K_i^{-1} + \frac{1}{2}  &- \frac{1}{2} & -1
    \end{bmatrix}.
    \end{align}
    \end{lemma}
     \begin{IEEEproof}
     Consider the incremental storage function
    \begin{align}\begin{split}
     Z_{2i}=& ~\frac{T_{\theta i}}{2}(\theta_i - \overline \theta_i^{opt})^2 + \frac{T_{si}}{2}(P_{si} - \overline P_{si}^{opt})^2 \nonumber
     \\&+\frac{T_{si}}{2}(P_{mi} - P_{si})^2 \\
    =& ~\frac{T_{\theta_i}}{2}(\theta_i - \overline \theta_i^{opt})^2 + T_{si}(P_{si} - \overline P_{si}^{opt})^2
     \\& +\frac{T_{si}}{2}(P_{mi} - \overline P_{mi}^{opt})^2 - T_{si}(P_{si} - \overline P_{si}^{opt})(P_{mi} - \overline P_{mi}^{opt}). 
\end{split}
  \end{align}
  It can be readily confirmed that $Z_{2i}$ is positive definite.
  We have that $Z_{2i}(P_{si}, \overline P_{si}, P_{mi}, \overline P_{mi}^{opt}, \theta_i, \overline \theta_i^{opt})$ satisfies along the solutions to (\ref{2tg}), (\ref{2tgc}),
  \begin{align}\label{dz2}
    \begin{split}
      \dot Z_{2i} =&~ (\theta_i - \overline \theta_i^{opt})(-\theta_i + P_{si}- (1 - K^{-1}_i)\omega_{gi})
      \\ &+ 2(P_{si} - \overline P_{si}^{opt})(-P_{si} - K^{-1}_i\omega_{gi}  + \theta_{i})
      \\ &+ T_{si}T_{mi}^{-1}(P_{mi} - \overline P_{mi}^{opt})(-P_{mi} + P_{si})
      \\ &- T_{si}T_{mi}^{-1}(P_{si} - \overline P_{si}^{opt})(-P_{mi} + P_{si})
      \\ &- (P_{mi} - \overline P_{mi}^{opt})( -P_{si} - K^{-1}_i\omega_{gi}  + \theta_{i})\\
      =& -T_{si}T_{mi}^{-1}(P_{si} - P_{mi})^2 - (P_{si} - \theta_i)^2 \\
      &-K_{i}^{-1}(P_{si} - P_{mi}) \omega_{gi} - K_{i}^{-1}(P_{si} - \theta_{i})\omega_{gi} \\
      &-(P_{si} - P_{mi})(P_{si} - \theta_i)\\
      &- (\theta_{i} - \overline \theta_i^{opt})\omega_{gi},
    \end{split}
  \end{align}
  where we exploited in the second identity the fact that at steady state
    \be\ba
   0  =& -\overline P_{si}^{opt} - K^{-1}_i 0+ \overline \theta_{i}^{opt} \\
    0   =& -\overline P_{mi}^{opt}  + \overline P_{si}^{opt} \\
  0 =& -\overline \theta_{i}^{opt}  + \overline P_{si}^{opt}  - (1 - K^{-1}_i)0,
  \ea
  \end{align}
  holds.
   We recall that $\dot U_{gi}=-D_{gi}\omega_{gi}^2 + \omega_{gi}(P_{mi} - \overline P_{mi}^{opt})$ and notice that
   \begin{align}
   \begin{split}
      &\quad \omega_{gi}(P_{mi} - \overline P_{mi}^{opt}) - (\theta_{i} - \overline \theta_i^{opt})\omega_{gi}\\
      =&\quad \omega_{gi}(P_{mi} - \theta_i) \\
      =&\quad \omega_{gi}(P_{si} - \theta_i)- \omega_{gi}(P_{si} - P_{mi}).
      \end{split}
   \end{align}
  The expression for $W_i$ then follows from writing $\dot U_{gi} + \dot Z_{2i}$ as a quadratic form.
     \end{IEEEproof}
     \bigskip
     We now address under what conditions $W_i$ is negative definite, which is important for the stability analysis in the next subsection.
     \begin{assumption}[Conditions on $K^{-1}_i$]\label{assumpdroop}
     Let the permanent droop constant $K_i$ be such that the following inequalities hold
     \be\label{assuma1}
     1 - \frac{T_{mi}}{T_{si}}  -\sqrt{\alpha_i} <  K_{i}^{-1} <  1 - \frac{T_{mi}}{T_{si}} + \sqrt{\alpha_i},
     \end{align}
     where
     \begin{align}\label{assuma2} \alpha_i = T_{mi}^2T_{si}^{-2}(4T_{si}T_{mi}^{-1} - 1)( D_{gi}T_{si}T_{mi}^{-1} - 1) . \end{align} Additonally, let $D_{gi}, T_{si}, T_{mi}$ be such that
     \begin{align}\label{assumineq}
     \begin{split}
     4T_{si}T_{mi}^{-1} &> 1\\
     D_{gi}T_{si}T_{mi}^{-1} &> 1,
     \end{split}
     \end{align}
   are satisfied.
     \end{assumption}
\begin{remark}[Locally verifiable]
  The power network generally consists of many generators. It is therefore important to note that the validity of Assumption \ref{assumpdroop} can be checked at each generator using only information that is locally available.
\end{remark}
    \begin{lemma}[Negative definiteness of $W_i$] Let Assumption 3 hold. Then $ W_i < \boldsymbol{0}$.
    \end{lemma}
    \begin{IEEEproof}
    Inequality (\ref{assumineq}) guarantees that \be X_i = \begin{bmatrix}
-T_{si}T_{mi}^{-1} &- \frac{1}{2} \\- \frac{1}{2} & -1
    \end{bmatrix} < \boldsymbol{0}.\end{align} It follows that $W_i <\boldsymbol{0}$ if and only if the Schur complement of $X_i$ in $W_i$ is negative definite. This Schur complement is given by
    \begin{align} \begin{split}
    S_i =& -D_{gi}  -\begin{bmatrix}
       -\frac{1}{2}K_i^{-1} -\frac{1}{2} \\ -\frac{1}{2}K_i^{-1} + \frac{1}{2}
    \end{bmatrix}^T X_i^{-1}\begin{bmatrix}
       -\frac{1}{2}K_i^{-1} -\frac{1}{2} \\  -\frac{1}{2}K_i^{-1} + \frac{1}{2}
    \end{bmatrix},
    \end{split}
     \end{align}
    and is quadratic in $K_i^{-1}$.  By Cramer's rule we have
       \begin{align}
         X_i^{-1} = \frac{1}{T_{si}T_{mi}^{-1} - \frac{1}{4}}\begin{bmatrix}
           -1 & \frac{1}{2} \\
           \frac{1}{2} & -T_{si}T_{mi}^{-1}
         \end{bmatrix},
       \end{align}
       and a straightforward calculation yields
        \begin{align}
       \begin{split}
     S_i =& -D_{gi} \\
     &+ \frac{\frac{1}{4}T_{si}T_{mi}^{-1}K_i^{-2}+(\frac{1}{2} -\frac{1}{2}T_{si}T_{mi}^{-1})K_i^{-1} + \frac{1}{2} + \frac{1}{4}T_{si}T_{mi}^{-1}}{T_{si}T_{mi}^{-1} - \frac{1}{4}}.
    \end{split}
    \end{align}
    The solution to $S_i=0$ is given by the quadratic formula resulting in
    \begin{align}
      K_i^{-1} =&~ \frac{-b_i}{2a_i} \pm \sqrt\frac{{b_i^2 - 4a_ic_i}}{4a^2_i},
    \end{align}
    with
    \begin{align}
    \begin{split}
      a_i =& ~ \frac{1}{4}T_{si}T_{mi}^{-1}\\
      b_i =& ~\frac{1}{2} -\frac{1}{2}T_{si}T_{mi}^{-1} \\
      c_i =& -D_{gi}(T_{si}T_{mi}^{-1} - \frac{1}{4}) + \frac{1}{2}+ \frac{1}{4}T_{si}T_{mi}^{-1}.
      \end{split}
    \end{align}
    Algebraic manipulations then yield
    \begin{align}\begin{split}
      \frac{-b_i}{2a_i} =& ~1 - \frac{T_{mi}}{T_{si}} \\
      \frac{b_i^2 - 4a_ic_i}{4a_i^2} =&~ T_{mi}^2T_{si}^{-2} - T_{mi}T_{si}^{-1}(4+D_{gi})+ 4D_{gi} \\
       =&~ T_{mi}^2T_{si}^{-2}(4T_{si}T_{mi}^{-1} - 1)( D_{gi}T_{si}T_{mi}^{-1} - 1) \\
      =& ~ \alpha_i.
      \end{split}
    \end{align}
       It can now be readily confirmed that $S_i < 0$ when (\ref{assuma1}) holds, where $\sqrt{\alpha_i}$ is real as a result of inequality (\ref{assumineq}).
    \end{IEEEproof}
    \bigskip

\subsection{Stability analysis and optimal distributed control}\label{subC}
Having discussed the separate control of the various turbine-governors, we now turn our attention to the question of how the different controllers in the network can cooperate to ensure minimization of the generation costs at steady state. To this end we add an additional communication term to controllers (\ref{1tgc}) and (\ref{2tgc}) \smod representing the exchange of information on the marginal costs among the controllers \emod
 \begin{align}
 \label{commcontol1}
 \begin{split}
    T_{\theta_{i}} \dot \theta_{i} =& -\theta_{i} + P_{mi}
     \\&  -K_{i}^{-1}q_i\sum_{\mathclap{j\in \mathcal{N}_{i}^{com}}} (q_i\theta_i + r_i- (q_j\theta_j + r_j)),
     \end{split} \quad \forall i \in \mathcal{V}_{g1}\\
     \begin{split}\label{commcontol2}
        T_{\theta_{i}} \dot \theta_{i} =& -\theta_{i} + P_{si} - (1 - K^{-1}_i)\omega_{gi} \\
    &  -q_i\sum_{\mathclap{j\in \mathcal{N}_{i}^{com}}} (q_i\theta_i + r_i- (q_j\theta_j + r_j)),
    \end{split} \quad \forall i \in \mathcal{V}_{g2}
  \end{align}
  where $\mathcal{N}_{i}^{com}$ is the set of buses connected via a communication link to bus $i$. The additional communication term can be interpreted as a consensus algorithm, where generator $i$ compares its marginal cost with the marginal costs of connected generators, such that the overall network converges to the state where there is consensus in the marginal costs (see Theorem 1). Due to the modified dynamics of the controller state $\theta_i$, the derivatives of $Z_{1i}$ and $Z_{21}$ along the solutions to (\ref{commcontol1}), (\ref{commcontol2}) need to be reevaluated. We exploit the result in the proof of Theorem 1, but is discussed separately for the sake of readability.
\begin{remark}[Communication induced modifications]\label{remarkLaplacian}
As a result of the additional communication term in (\ref{commcontol1}), (\ref{commcontol2}), the expressions for $\dot Z_{1i}$ and $\dot Z_{2i}$ given in respectively (\ref{zi}) and (\ref{dz2}) need to be modified.
  Notice that \begin{align} q_i\sum_{\mathclap{j\in \mathcal{N}_{i}^{com}}} (q_i\theta_i + r_i- (q_j\theta_j + r_j))= \big(QL^{com}(Q\theta + R)\big)_i,
   \end{align} where $L^{com}$ is the Laplacian matrix reflecting the topology of the communication network.
   Therefore, we add the following term to $\dot Z_{1i}$ and $\dot Z_{2i}$
   \begin{align}
   -(\theta_i - \overline \theta_i^{opt})\big(QL^{com}(Q\theta + R)\big)_i
   \end{align}
   Summing over all buses $i \in \mathcal{V}_g$ then yields
    \begin{align}
    \begin{split}
    &-\sum_{i \in \mathcal{V}_g}(\theta_i - \overline \theta_i^{opt})\big(QL^{com}(Q\theta + R)\big)_i\\
   =&-(\theta- \overline \theta^{opt})^TQL^{com}(Q\theta + R) \\
   =& -(Q\theta + R - (Q\overline \theta^{opt} + R))^TL^{com}\\&\quad \cdot(Q\theta + R - (Q\overline \theta^{opt} + R)),
   \end{split}
   \end{align}
   where we exploited
   \begin{align}
     L^{com}(Q\overline \theta^{opt} + R) = \boldsymbol{0},
   \end{align}
   which is a result of $\overline \theta^{opt} = \overline P_m^{opt}$, $Q\overline \theta^{opt} + R \in {\rm Im}(\mathds{1}_{n_g})$ and ${\rm Ker}(L^{com}) = {\rm Im}(\mathds{1}_{n_g})$.
\end{remark}

 The communication network is utilized to ensure that all marginal costs converge to the same value throughout the network (see the proof of Theorem 1), leading to the following assumption:
\begin{assumption}[Connectivity]\label{assum3}
  The graph reflecting the topology of information exchange among the controllers is undirected and connected, but can differ from the topology of the power network.
\end{assumption}
We are now ready to state the main result of this work.
\begin{theorem}[Distributed optimal LFC]
  Let assumptions 1, 2, 3 and 4 hold. Consider the power network (\ref{syscompact}), turbine-governor dynamics (\ref{1tg}), (\ref{2tg}) and the distributed controllers (\ref{commcontol1}), (\ref{commcontol2}). Then, solutions that start sufficiently close to $(\overline \eta, \overline \omega = \boldsymbol{0}, \overline P_m^{opt}, \overline P_s^{opt}, \overline \theta^{opt})$ converge to the set where we have frequency regulation and where the power generation solves optimization problem (\ref{optimal}), i.e. $\overline \omega = \boldsymbol{0}$ and $\overline P_m = \overline P_m^{opt}$.
  \end{theorem}
\begin{IEEEproof}
  As a result of Lemma \ref{u1}, Corollary 1, Lemma \ref{lemma2ndorder} and Remark \ref{remarkLaplacian}, we have that $U + \sum_{i \in \mathcal{V}_{g1}}Z_{1i} + \sum_{i \in \mathcal{V}_{g2}}Z_{2i}$ satisfies
  \begin{align}\label{totalUdot}
   \begin{split}
   &\dot U + \sum_{i \in \mathcal{V}_{g1}}\dot Z_{1i} + \sum_{i \in \mathcal{V}_{g2}}\dot Z_{2i}\\
    &= -\|\omega_l \|^2_{D_l} + \sum_{i \in \mathcal{V}_{g1}}\Big( -D_{gi}\omega_{gi}^2  -K_i(\theta_i - P_{mi})^2 \Big) \\
   &+  \sum_{i \in \mathcal{V}_{g2}}  \begin{bmatrix}
      \omega_{gi} \\ P_{si} - P_{mi} \\ P_{si} - \theta_i
    \end{bmatrix}^T  W_i \begin{bmatrix}
      \omega_{gi} \\ P_{si} - P_{mi} \\ P_{si} - \theta_i
    \end{bmatrix} \\
    &-(Q\theta + R - (Q\overline \theta^{opt} + R))^TL^{com}(Q\theta + R - (Q\overline \theta^{opt} + R))\\&\leq 0,
    \end{split}
  \end{align}
  along the solutions to the power network (\ref{syscompact}), turbine-governor dynamics (\ref{1tg}), (\ref{2tg}) and the distributed controllers (\ref{commcontol1}), (\ref{commcontol2}).
  Particularly, it follows from Assumption 3 that $W_i < 0$. Since $(\overline \eta, \overline \omega = \boldsymbol{0}, \overline P_m^{opt}, \overline P_s^{opt}, \overline \theta^{opt})$ is a strict local minimum of $U + \sum_{i \in \mathcal{V}_{g1}}Z_{1i} + \sum_{i \in \mathcal{V}_{g2}}Z_{2i}$ as a consequence of Assumption 2, there exists a compact level set
 $
  \Upsilon$ around  $(\overline \eta, \overline \omega = \boldsymbol{0}, \overline P_m^{opt}, \overline P_s^{opt}, \overline \theta^{opt})$, which is forward invariant. By LaSalle's invariance principle, any solution starting in $\Upsilon$ asymptotically converges to the largest invariant set contained in \begin{align}\begin{split}
  &\Upsilon \cap \{(\eta, \omega,  P_m, P_s, \theta): \\ & \omega = \boldsymbol{0}, ~P_m = \theta, ~Q\theta + R = Q\overline \theta^{opt} + R + c \mathds{1} \},\end{split}\end{align}
  where $c \in \mathds{R}$ is a scalar, and $Q\theta + R = Q\overline \theta^{opt} + R + c\mathds{1}$ follows from the connectedness of the communication graph. Since $P_m = \theta = \overline \theta^{opt} +cQ^{-1}\mathds{1} = \overline P_m^{opt} +cQ^{-1}\mathds{1}$, the power network satisfies on this invariant set
  \be\ba\label{syscompactsteadylasalle}
   \dot \eta =&~\mathcal{B}^T \boldsymbol{0} \\
   \boldsymbol{0}  =& -D_{g}\boldsymbol{0}  -  \mathcal{B}_{g} \Gamma \sin( \eta) + \overline P_m^{opt} + c Q^{-1}\mathds{1} \\
   \boldsymbol{0} =&  -D_{l}\boldsymbol{0}  -  \mathcal{B}_{l} \Gamma \sin(\eta) - P_l. \\
  \ea
  \end{align}
  Premultiplying the second and third line of (\ref{syscompactsteadylasalle}) with $\mathds{1}_n^T$, we have
    \be \mathds{1}_n^T \begin{bmatrix}
      -D_{g}\boldsymbol{0}  -  \mathcal{B}_{g} \Gamma \sin( \eta) + \overline P_m^{opt} + cQ^{-1}\mathds{1} \\
     -D_{l}\boldsymbol{0}  -  \mathcal{B}_{l} \Gamma \sin(\eta) - P_l.  \end{bmatrix} = 0.
   \end{align}
  Since $\mathds{1}_n^T \begin{bmatrix}
    \mathcal{B}_g \\ \mathcal{B}_l
  \end{bmatrix} = \boldsymbol{0}$, $\mathds{1}_{n_g}^T \overline P^{opt}_m - \mathds{1}_{n_l}^T P_l = 0$ and $Q^{-1}$ is a diagonal matrix with only positive elements, it follows that necessarily $c = 0$ and therefore $\theta = \overline \theta^{opt}$. We can conclude that the system indeed converges to the set where $\omega = \boldsymbol{0}$ and $P_m = \overline P_m^{opt}$, characterized in Lemma \ref{lemma4}.
\end{IEEEproof}
\smod
\begin{remark}[Region of attraction]
The local nature of our result is a consequence of the  considered incremental storage function having  a local minimum at the desired steady state. Nevertheless, the provided results are helpful to further characterize various sublevel sets of the incremental storage function (\cite{dvijotham_2015_arxiv, vu_2016_tps,persis_2016_arxiv_microgrid}), for instance by numerically assessing the sublevel sets that are compact. We leave
   a thorough analysis of the region of attraction as an interesting future direction.
\end{remark}

\emod
\begin{remark}[Primal-dual based approaches]
A popular alternative to the consensus based algorithm (\ref{commcontol1}), (\ref{commcontol2}) is a primal-dual gradient based approach.
To obtain a distributed solution, optimization problem (\ref{optimal}) is replaced\footnote{See \cite[Lemma 4]{trip_2016_automatica} for a discussion on the equivalence of (\ref{optimal}) and (\ref{optimalpd}).} by
\be\ba\label{optimalpd}
& \min_{P_m} C(P_m)\\
{\rm s.t.}\quad & \boldsymbol{0}= -\mathcal{B}v +\begin{bmatrix}
  P_m \\ -P_l
\end{bmatrix}.
\ea
\end{align}
The associated Lagrangian function is given by
\begin{align}
L(P_m, \lambda) = C(P_m) + \lambda^T\Bigg(-\mathcal{B}v +\begin{bmatrix}
  P_m \\ -P_l\end{bmatrix}\Bigg),
\end{align}
where $\lambda$ is called the Lagrange multiplier. Under convexity of (\ref{optimalpd}), strong duality holds and the solution to (\ref{optimalpd}) is equivalent \cite{boyd_2004_book} to the solution to
\begin{align}\label{dual}
  \max_\lambda \min_{P_m} L(P_m, \lambda).
\end{align}
Following \cite{zhang_2015_automatica, chen2015, stegink_2016_arxiv}, a continuous primal-dual algorithm can be exploited to solve (\ref{dual}).
However, since the evolution of $P_m$ is described by the turbine dynamics, we cannot design its dynamics.
Bearing in mind that controller (\ref{1tgc}) and (\ref{2tgc}) enforce a steady state where $P_m = \theta$, we solve instead
\begin{align}
  \max_\lambda \min_{\theta} L(\theta, \lambda),
\end{align}
where the dynamics of $\theta$ can be freely adjusted.
Inspired by the results in \cite{zhang_2015_automatica, chen2015, stegink_2016_arxiv}, we replace the communication term in (\ref{commcontol1}), (\ref{commcontol2}),
\begin{align}-q_i\sum_{\mathclap{j\in \mathcal{N}_{i}^{com}}} (q_i\theta_i + r_i- (q_j\theta_j + r_j))\end{align}
by
  \begin{align}
    \frac{\partial L}{\partial \theta_{i}} = -\nabla C_i(\theta_i) + \lambda_i,
  \end{align}
yielding the modified controllers
 \begin{align}
 \begin{split}
    T_{\theta_{i}} \dot \theta_{i} =& -\theta_{i} + P_{mi}
      -K_{i}^{-1}(\nabla C_i(\theta_i) - \lambda_i) ,
     \end{split} ~~ \forall i \in \mathcal{V}_{g1}\\
     \begin{split}
        T_{\theta_{i}} \dot \theta_{i} =& -\theta_{i} + P_{si} - (1 - K^{-1}_i)\omega_{gi}
      \\&-(\nabla C_i(\theta_i) - \lambda_i).
    \end{split} \quad ~~~~~~~~\forall i \in \mathcal{V}_{g2}
  \end{align}
  The variables $v$ and $\lambda$ evolve according to
  \begin{align}\label{odemultiplier}
  \begin{split}
  \dot v =&~ \frac{\partial L}{\partial v} = -\mathcal{B}^T\lambda \\
  \dot \lambda =& -\frac{\partial L}{\partial \lambda} = \mathcal{B}v -\begin{bmatrix}
  \theta \\ -P_l\end{bmatrix}.
  \end{split}
  \end{align}
  The analysis of Theorem 1 can now be repeated with the additional storage term
  \begin{align}
    Z_3 =& ~\frac{1}{2}(v - \overline v)^T(v - \overline v)  +  \frac{1}{2}(\lambda - \overline \lambda)^T(\lambda - \overline \lambda).
  \end{align}
  We notice that in this case only convexity of $C(\cdot)$ is required and that the load $P_l$ appears in (\ref{odemultiplier}).
\end{remark}
\begin{remark}[Load control]\label{remarkload}
Incorporating load control in the LFC has been recently studied in e.g. \cite{zhao_2013_ps, chen_2012_book, weckx_2015_tps} and can be incorporated within the presented framework with minor modifications with respect to the previous discussion.
To do so, we modify the dynamics at the load buses $i \in \mathcal{V}_{l}$  to become
  \be\ba\label{newload}
  \dot{\delta}_{i} =&~\omega_{li}  \\
   0 =& -D_{li}\omega_{li} \\& -  \sum_{j \in \mathcal{N}_{i}} V_{i}V_{j}B_{ij} \sin(\delta_{i} - \delta_{j}) - P_{li} - u_{li},
  \ea
  \end{align}
where $u_{li}$ is the additional controllable load. Associated to every controllable load is a strictly concave benefit function of the form \be
C^{B}_i(u_{li}) = \frac{1}{2}q_{i}u_{li}^2 + r_{i}u_{li} + s_{i},
\end{align}
which is a common approach to quantify the benefit of the consumed power. Instead of minimizing the total generation costs as in (\ref{optimal}) we now aim at maximizing the so-called `social welfare' \cite{berger_1989_tps}, \cite{kiani_2011_cdc},
\be\ba\label{socialwelfare}
& \max_{u_l, P_m} C^B(u_l) -  C(P_m)\\
{\rm s.t.} \quad& 0= \mathds{1}_{n_g}^T \overline P_{m}- \mathds{1}_{n_l}^T (P_{l} + u_l),
\ea
\end{align}
where $C^B(u_l) -  C(P_m) = \sum_{i \in \mathcal{V}_{l}} C^B_i(u_{li}) - \sum_{i \in \mathcal{V}_{g}} C_i(P_{mi})$. Notice that (\ref{socialwelfare}) is equivalent to (\ref{optimal}) in the absence of controllable loads. A straightforward but remarkable extension of Lemma \ref{u1} is that $U(\eta, \overline\eta, \omega_{g}, \overline \omega_{g})$ as in (\ref{storagefunction}) now satisfies along the solutions to (\ref{syscompact}) and (\ref{newload})
\begin{align}\begin{split}
  \dot U =& -\|\omega_{g} - \overline \omega_{g}\|_{ D_{g}}^2  -\|\omega_{l} - \overline \omega_{l}\|_{ D_{l}}^2 \\& + (\omega_{g} - \overline \omega_{g})^T(P_m - \overline P_m) \\&- (\omega_{l} - \overline \omega_{l})^T(u_l - \overline u_l),
\end{split}\end{align}
i.e. the power network is also output strictly cyclo-incrementally passive with respect to the additional input-output pair $(u_l, -\omega_l)$.
This property allows to incorporate load control in the same manner as the generation control.
A thorough discussion on all possible load dynamics is outside the scope of this paper, although the considered turbine-governor dynamics can be straightforwardly adapted.
   In the case there are no restrictions on the design, a possible load controller is given by
    \begin{align}
    \begin{split}
        T_{\theta_{i}} \dot \theta_{i} =&~ \omega_{li} - q_i\sum_{j\in \mathcal{N}_{i}^{com}} (q_i\theta_i + r_i- (q_j\theta_j + r_j)) \\
        u_{li} = ~&\theta_{li}.
    \end{split}  \quad \forall i \in \mathcal{V}_{l} \end{align}
     The analysis of Theorem 1 can now be repeated with the additional storage term
 \begin{align}
 Z_{3i} = \frac{1}{2}\sum_{i \in \mathcal{V}_l}(\theta_i - \overline \theta_i^{opt})^2.
 \end{align}
\end{remark}
\smod
\begin{remark}[Time-varying loads]
  Theorem 1 above establishes frequency regulation under the assumption of a constant unknown load $P_l$. In a realistic setting the (net) load, including uncontrollable renewable energy generation, is likely to change erratically. Although exact frequency regulation is not possible in that case, the results in this paper are useful to bound the resulting frequency deviation.  If a varying load $Q_l(t)$ with finite $\mathcal{L}_2$-norm ($\int_0^{\infty}\|Q_l(\tau) \|^2 d\tau < \infty$) is added to the load bus, e.g. by taking $u_{li} = Q_{li}(t)$ in (\ref{newload}), it is possible, following \cite[Remark 8]{trip_2016_automatica}, to derive from (\ref{totalUdot}) the existence of a finite $\mathcal{L}_2$-to-$\mathcal{L}_\infty$ gain and a finite $\mathcal{L}_2$-to-$\mathcal{L}_2$ gain from the load (disturbance) $Q_l$ to the frequency deviation $\omega$ \cite{kundur2004}.
\end{remark}
\emod

\section{Case study}\label{sec7}
To illustrate the proposed control scheme we adopt the 6 bus system from \cite{wood1996power}. Its topology is shown in Figure 1. The relevant generator and load parameters are provided in Table 2, whereas the transmission line parameters are provided in Table 3. The used numerical values are based on \cite{wood1996power} and \cite{venkat_2008_cst}. The turbine-governor dynamics are modelled by the second-order model (\ref{2tg}). Every generator is equipped with the controller presented in (\ref{commcontol2}). The communication links between the controllers are also depicted in Figure 1. The system is initially at steady state with loads $P_{l1}, P_{l2}$ and $P_{l3}$ being 1.01, 1.20 and 1.18 pu respectively (assuming a base power of 100 MVA). After 10 seconds the loads are respectively increased to $1.15$, 1.25 and 1.21 pu. From Figure 2 we can see how the controllers regulate the frequency deviation back to zero. The total generation is shared optimally among the different generators such that (\ref{optimal}) is solved.

\newcommand{\genset}[3]{
  \node[circle,draw,thick,minimum width=7mm,inner sep=0pt,drop
  shadow,fill=white,#3] (#1) at (#2) {};
  \draw[thick] ($(#2)-(2mm,0)$) sin ++(1mm,1mm) cos ++(1mm,-1mm) sin
  ++(1mm,-1mm) cos ++(1mm,1mm); }

\begin{figure}
\medskip
\centering
\begin{tikzpicture}[thick, scale=0.8, transform shape]

\tikzstyle{line}=[-, thick]
\tikzstyle{loadline}=[->,thick,>=stealth']
\tikzstyle{busbar} = [rectangle,draw,fill=black,inner sep=0pt];
\tikzstyle{hbus} = [busbar,minimum width=10mm,minimum height=2pt,drop shadow];

\coordinate (c1) at (0,3);
\coordinate (c2) at (3.5,4.5);
\coordinate (c3) at (7,3);
\coordinate (c4) at (0,0);
\coordinate (c5) at (3.5,-1.5);
\coordinate (c6) at (7,0);
\coordinate (over) at (0,3.75);

\begin{footnotesize}
\node[hbus,minimum width=16mm,label=left:Bus 1] (b1) at (c1) {};
\node[hbus,minimum width=25mm,label=left:Bus 2] (b2) at (c2) {};
\node[hbus,minimum width=16mm,label=right:Bus 3] (b3) at (c3) {};
\node[hbus,minimum width=16mm,label=left:Bus 4] (b4) at (c4) {};
\node[hbus,minimum width=25mm,label=left:Bus 5] (b5) at (c5) {};
\node[hbus,minimum width=16mm,label=right:Bus 6] (b6) at (c6) {};
\end{footnotesize}

\draw[line,<->, dashed, blue] ([xshift=-3mm, yshift=11mm] b1.north) to[bend left] ([xshift=-3mm, yshift=11mm] b2.south);
\draw[line,<->, dashed, blue] ([xshift=3mm, yshift=11mm] b3.north) to[bend right] ([xshift=3mm, yshift=11mm] b2.south);
\draw[line] ([xshift=5mm] b1.north) |- ([xshift=-10mm,yshift=-5mm] b2.south)  --
([xshift=-10mm] b2.south);
\draw[line] ([xshift=-5mm] b1.south) -- ([xshift=-5mm] b4.north) ;
\draw[line] ([xshift=5mm] b1.south) -- ([xshift=5mm,yshift=-5mm] b1.south)
-- ([xshift=-5mm,yshift=8mm] b5.north) -- ([xshift=-5mm] b5.north);
\draw[line] ([xshift=10mm] b2.south) -- ++(0,-0.5) -| ([xshift=-5mm]
b3.north);
\draw[line] ([xshift=-5mm] b2.south) --
([xshift=-5mm,yshift=-8mm] b2.south) -- ([xshift=5mm,yshift=5mm] b4.north)
-- ([xshift=5mm] b4.north);
\draw[line] (b2.south) -- (b5.north);
\draw[line] ([xshift=5mm] b2.south) -- ([xshift=5mm,yshift=-8mm] b2.south)
-- ([xshift=-5mm,yshift=5mm] b6.north) -- ([xshift=-5mm] b6.north);
\draw[line] ([xshift=-5mm] b3.south) -- ([xshift=-5mm,yshift=-5mm] b3.south)
-- ([xshift=5mm,yshift=8mm] b5.north) -- ([xshift=5mm] b5.north);
\draw[line] ([xshift=5mm] b3.south) -- ([xshift=5mm] b6.north);
\draw[line] ([xshift=5mm] b4.south) |- ([xshift=-10mm,yshift=5mm] b5.north) --
([xshift=-10mm] b5.north);
\draw[line] ([xshift=10mm] b5.north) -- ([xshift=10mm,yshift=5mm] b5.north) -|
([xshift=-5mm] b6.south);

\genset{g1}{$(c1)+(-5mm,8mm)$}{label=above:$g_1$}
\draw[line] ([xshift=-5mm] b1.north) -- (g1.south);
\genset{g2}{$(c2)+(0,8mm)$}{label=above:$g_2$}
\draw[line] (b2.north) -- (g2.south);
\genset{g3}{$(c3)+(5mm,8mm)$}{label=above:$g_3$}
\draw[line] ([xshift=5mm] b3.north) -- (g3.south);

\draw[loadline] ([xshift=-5mm] b4.south) -- ++(0,-0.8) node[text centered,text
width=10mm,below] {$l_4$};
\draw[loadline] (b5.south) -- ++(0,-0.8) node[text centered,text
width=10mm,below] {$l_5$};
\draw[loadline] ([xshift=5mm] b6.south) -- ++(0,-0.8) node[text centered,text
width=10mm,below] {$l_6$};

\end{tikzpicture}
\caption{Diagram for a 6 bus power network, consisting of 3 generator and 3 load buses. The turbine-governor dynamics of generators are represented by a second-order model. The communication links are represented by the dashed lines.}
\label{fig:case6ww}
\end{figure}
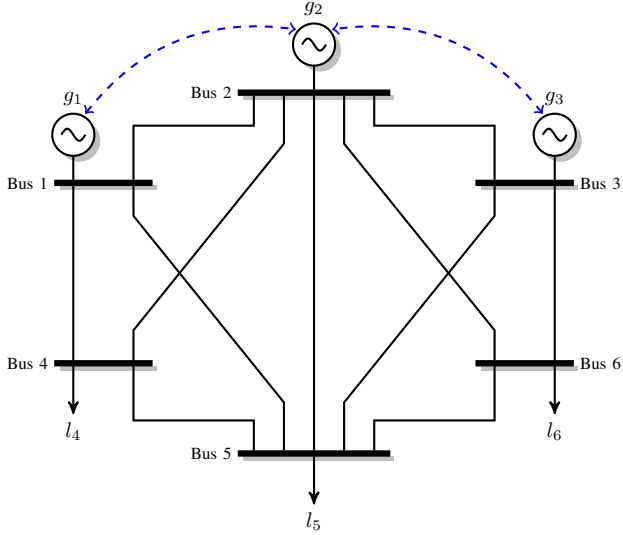

\begin{table}\label{tab2} \centering
\scalebox{1}{\tabcolsep=0.13cm
    \begin{tabular}{l  |c c c c c c }
         & \rot{Bus 1} & \rot{Bus 2} & \rot{Bus 3} & \rot{Bus 4} & \rot{Bus 5} & \rot{Bus 6}\\
        \hline
 $M_{i}$ & 4.6& 6.2& 5.1 & -- & -- & --\\
  $D_{gi}$  &  3.4& 3.0 & 4.2&-- & -- & --\\
  $D_{li}$   & --& -- & --& 1.0 & 1.6 & 1.2\\
  $V_{i}$ &  1.05& 0.98 & 1.04& 1.01 & 1.03 & 1.00\\
  $T_{si}$   & 4.0& 4.6& 5.0 & -- & -- & --\\
  $T_{mi}$   & 5.0& 6.7& 10.0 & -- & -- & --\\
  $K_{i}$   & 0.5& 0.5 & 0.5 & -- & -- & --\\
  $T_{\theta i}$   & 0.1& 0.1& 0.1& -- & -- & --\\
  $q_{i}$   & 2.4& 3.8 & 3.4 & -- & -- & --\\
  $r_i$  & 10.5& 5.7 & 8.9& -- & -- & --\\
  $s_i$  & 9.1& 14.4 & 13.2 & -- & -- & --\\
    \end{tabular}}
    \caption*{Table 2: Numerical values of the generator and load parameters. The values for $K_i$ satisfy Assumption 3.}
\end{table}
\begin{table}\label{tab2} \centering
    \scalebox{1}{\begin{tabular}{c| c c c c c c c}
        $B_{ij}$ (pu)&  1 & 2 & 3 & 4 & 5 & 6 & $j$\\
     \hline
  1 & --& -4.0& -- & -4.7 & -3.1 & --\\
  2& -4.0& -- & -3.8& -8.0 & -3.0 & -4.5\\
  3&  --& -3.8 & --& -- & -3.2 & -9.6\\
  4 & -4.7& -8.0& -- & -- & -2.0 & --\\
  5& -3.1& -3.0 & -3.2 & -2.0 & -- & -3.0\\
  6& --& -4.5&-9.6 & -- & -3.0 & --\\
  $i$
    \end{tabular}}
    \caption*{Table 3: Susceptance $B_{ij}$ of the transmission line connecting bus $i$ and bus $j$. Values are per unit on a base of 100 MVA.}
\end{table}
 \begin{figure}\label{figth5}
\centering
  \includegraphics[trim= 2.4cm 7cm 2.4cm 7cm, width=\columnwidth]{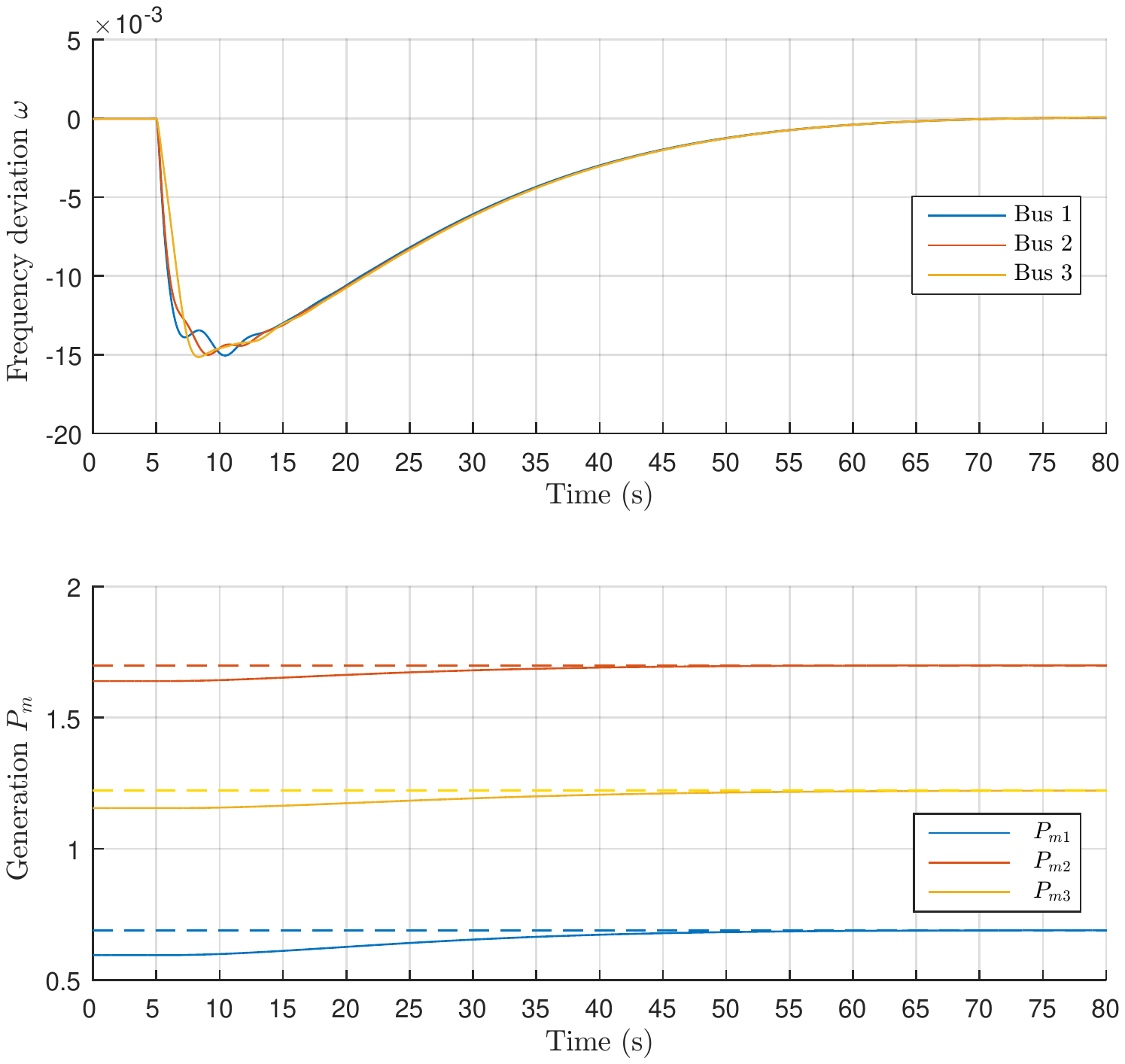}
  \caption{Frequency response and generated power at the generator buses using the controllers (\ref{commcontol2}). The load is increased at timestep 5, whereafter the frequency deviation is regulated back to zero and generation costs are minimized.  The cost minimizing generation $\overline P_m^{opt}$ for $t > 5$, characterized in Lemma \ref{lemma4}, is given by the dashed lines.
  }
  \end{figure}

  \subsection{Instability}
We now show that a wrongly chosen value for the frequency gain $(1 - K_i^{-1})$ in controller (\ref{commcontol2}) can lead to instability. To do so, we change the controller at generator $3$ into
  \begin{align}
 \begin{split}\label{unstable}
        T_{\theta_{3}} \dot \theta_{3} =& -\theta_{3} + P_{s3} - 5(1 - K^{-1}_3)\omega_{g3} \\
    &  -q_3\sum_{\mathclap{j\in \mathcal{N}_{3}^{com}}} (q_3\theta_3 + r_3- (q_j\theta_j + r_j)),
    \end{split}
  \end{align}
 for $t> 5$.
Leaving all other values identical to the previous simulation, we notice from Figure 3 that this change at only one generator can cause instability throughout the whole network.
 \begin{figure}\label{figth5}
\centering
  \includegraphics[trim= 2.4cm 7cm 2.4cm 7cm, width=\columnwidth]{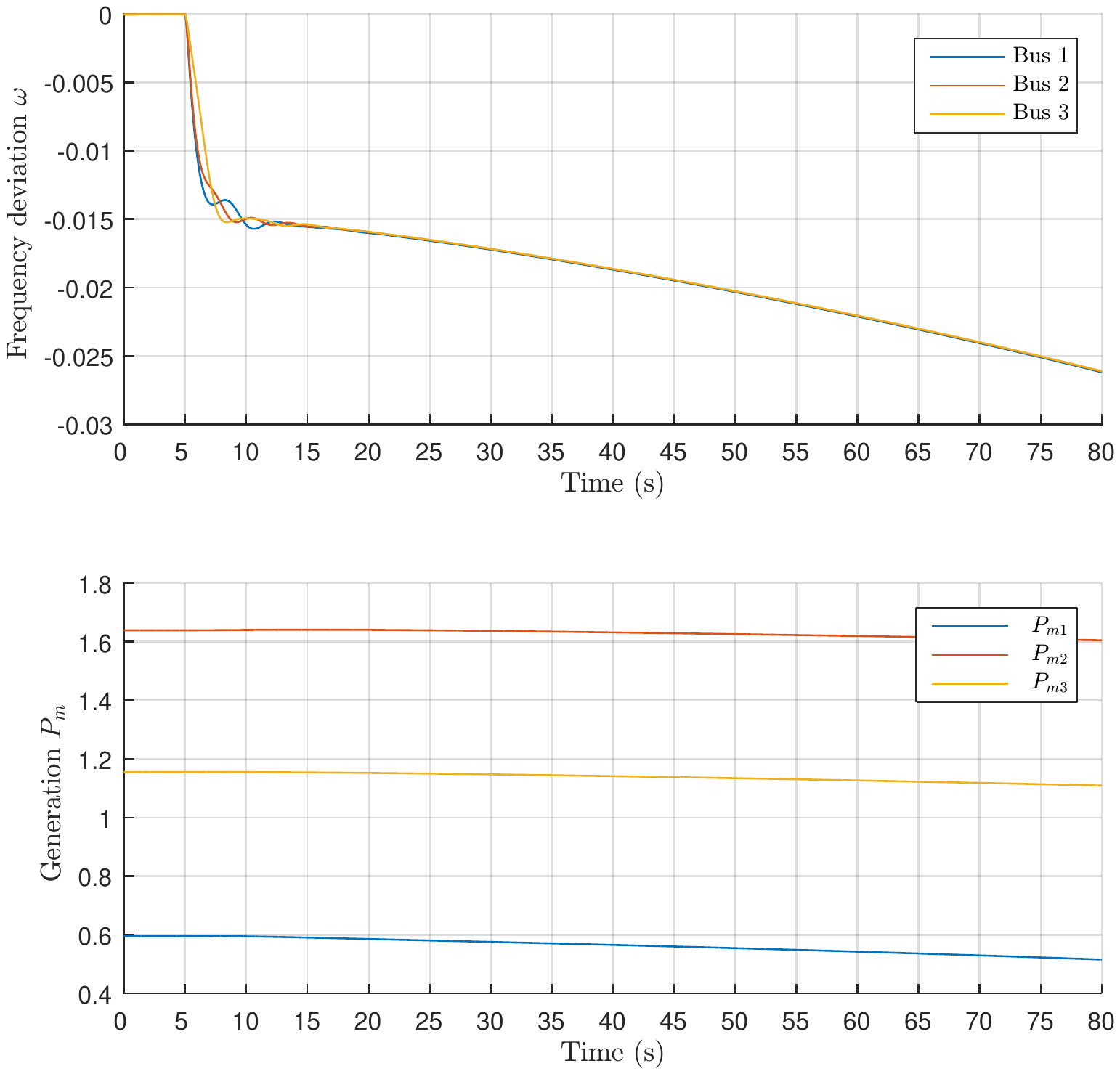}
  \caption{Frequency response and generated power at the generator buses using the controllers (\ref{commcontol2}). The load is increased at timestep 5 and the controller at generator 3 is replaced by (\ref{unstable}). Both the frequency deviation and the power generation become unstable.
  }
  \end{figure}

\section{Conclusions and future research}\label{sec8}
We presented the design of a distributed optimal LFC control architecture that regulates the frequency in the power network, while minimizing the generation costs (or maximizing the social welfare in case of controllable loads). Based on an energy function of the power network we derived an incremental passivity property for a well known structure preserving network model. The passivity property then facilitates the design of distributed controllers that adjust the input to the turbine-governor and load. In this work we have considered a first-order and a (non-passive) second-order model describing the turbine-governor dynamics. We establish a locally verifiable range of acceptable droop constants for the second-order model that allows us to infer stability. The presented results provide the opportunity to include the important turbine-governor models in the stability analysis of optimal LFC in a more realistic manner than was previously possible.
\medskip \\
There are various interesting extensions to the presented work. We briefly discuss a few.
\smod
The  distributed controllers (\ref{commcontol1}), (\ref{commcontol2}) are shown to solve optimization problem (\ref{optimal}). Minimizing general convex cost functions and satisfying transmission and generator constraints using consensus based controllers is still an open problem.
\emod
The distributed control architecture employs a communication network to obtain the desired optimality features. In this work the communication is assumed to be continuous and instantaneous. It is desirable to relax these communication assumptions. Lyapunov arguments have been used to design distributed event-triggered control algorithms in e.g. \cite{persis_2016_TAC}, \cite{postoyan2015} within a hybrid system framework. A promising research direction is to adapt these results to the present setting. \smod Furthermore, we note that Assumption 3 provides a sufficient condition on the parameters of the second-order turbine-governor model to infer stability of the overall network. Exploring the necessity of this condition, and potentially relaxing it, might offer additional insights on the role of turbine-governor dynamics within (optimal) LFC schemes.  \emod

\balance
\bibliographystyle{IEEEtran}
\bibliography{bibliography}
\end{document}